# OPTIMAL RATE OF CONVERGENCE FOR NONPARAMETRIC CHANGE-POINT ESTIMATORS FOR NONSTATIONARY SEQUENCES

By Samir Ben Hariz, Jonathan J. Wylie and Qiang Zhang[1]

*Université du Maine, City University of Hong Kong
and City University of Hong Kong*

Let $(X_i)_{i=1,\ldots,n}$ be a possibly nonstationary sequence such that $\mathscr{L}(X_i) = P_n$ if $i \leq n\theta$ and $\mathscr{L}(X_i) = Q_n$ if $i > n\theta$, where $0 < \theta < 1$ is the location of the change-point to be estimated. We construct a class of estimators based on the empirical measures and a seminorm on the space of measures defined through a family of functions $\mathcal{F}$. We prove the consistency of the estimator and give rates of convergence under very general conditions. In particular, the $1/n$ rate is achieved for a wide class of processes including long-range dependent sequences and even nonstationary ones. The approach unifies, generalizes and improves on the existing results for both parametric and nonparametric change-point estimation, applied to independent, short-range dependent and as well long-range dependent sequences.

**1. Introduction.** The change-point problem, in which one must detect a change in the marginal distribution of a random sequence, is important in a wide range of applications and has therefore become a classical problem in statistics. A comprehensive review of the subject can be found in [5]. In this paper we consider the general case of nonparametric estimation that must be used when no a priori information regarding the marginal distributions before and after the change-point is known. Although this problem has been widely studied for independent sequences, studying dependent sequences has importance for both theoretical reasons and numerous practical applications. In this paper we consider this challenging problem and develop a unified framework in which we can deal with sequences with quite general dependence structures. We prove that the rate of convergence of a

Received December 2005; revised September 2006.
[1]Research supported by City University of Hong Kong contract 7001561.

broad family of nonparametric estimators is $O_p(n^{-1})$. This is a particularly surprising result because the dependence structure of the sequence plays absolutely no role in determining the rate of convergence. The rate $O_p(n^{-1})$ is clearly optimal because there are only $n$ points in the sequence.

For independent sequences there is a wide literature, and both parametric and nonparametric methods have been widely studied. The nonparametric problem was considered by Carlstein [4], who proposed an estimator, proved its consistency and determined a rate of convergence. Dümbgen [6] embedded the estimator proposed by Carlstein in a more general framework, improved the rate of convergence in probability and derived the limiting distributions for certain models. Ferger [7] considered the almost sure convergence for Dümbgen's estimators. Yao, Huang and Davis [15] considered the case in which the location of the change-point can tend to either 0 or 1 as the sequence length tends to infinity. Ferger [8, 9] has investigated a number of features of change-point estimators including probability bounds and rates of weak and almost sure convergence. Since then several works have generalized these results to a weakly dependent or short-range dependent setting.

In recent years the importance of long-memory or long-range dependent (LRD) processes has been realized in a wide range of applications, especially in the analysis of financial and telecommunication data. For the purposes of this paper we define real sequences $(X_i)_{i=1,\ldots,n}$ to be short-range dependent (SRD) if $\limsup_{n\to\infty} n^{-1} \mathbb{E}[\sum_{i=1}^{n}(X_i - \mathbb{E}[X_i])]^2 < \infty$ and LRD otherwise. Several works are concerned with the generalization of the results for independent sequences to a SRD setting. However, estimating change-points for LRD sequences poses a number of significant challenges and there are much fewer known results in this case.

Parametric change-point estimation for LRD sequences, in which one typically has a priori knowledge about the marginal distributions, has been considered by a number of authors. Kokoszka and Leipus [12] considered the change in the mean for dependent observations for LRD sequences. They obtained rates in probability for the cumulative sum (CUSUM) change-point estimator and gave a rate of convergence of the estimator that gets worse as the strength of the dependence increases. The problem with a jump in the mean that tends to zero was considered by Horváth and Kokoszka [11]. They proved the consistency of the estimator and gave the limiting distribution. For sequences that have a change in the mean, Ben Hariz and Wylie [2] showed that the rate of convergence does not get worse as the strength of the dependence increases and that the rate of convergence for independent sequences is also achieved for both SRD and LRD sequences. In the nonparametric setting Giraitis, Leipus and Surgailis [10] derived a number of results that focused mainly on hypothesis testing. However, to our knowledge, there are no results regarding rates of convergence of nonparametric change-point estimation for LRD sequences.

In this paper we adopt a very general framework that allows us to consider a very general class of dependence structures. In particular, we make no assumption about stationarity in the dependence structure. This is especially important in practice because one can confidently make use of the proposed estimators on a sequence without checking for such stationarity (which is typically extremely difficult in practice). This framework represents a unified setting in which independent, SRD and LRD sequences can be treated. We prove the consistency of a Dümbgen-type estimator and show that the $O_p(n^{-1})$ rate of convergence for independent sequences is also achieved for both SRD and LRD sequences. In addition, we consider the case in which the difference between the distributions before and after the change-point tends to zero.

**2. Main results.** Let $(X_i)_{i=1,\ldots,n}$ be a sequence in a measurable space $E$. The marginal distribution (which may depend on the sequence length $n$) is given by

$$\mathscr{L}(X_i) = \begin{cases} P_n, & \text{if } i \leq n\theta, \\ Q_n, & \text{if } i > n\theta, \end{cases}$$

where $0 < \theta < 1$ is the location of the change-point. This means that we assume first-order stationarity on either side of the change-point, but make no assumption about stationarity in the dependence structure of the sequence.

Given the sequence $(X_i)_{i=1,\ldots,n}$, we aim to estimate the location of the change-point $\theta$ using an estimator of the general type

$$(2.1) \qquad \hat{\theta}_n = \frac{1}{n} \min\left(\arg\max_{1 \leq k < n}\{N(D_k)\}\right),$$

where $N$ is a (possibly random) seminorm on the space $\mathcal{M}$ of signed finite measures on $E$,

$$(2.2) \qquad D_k = \left[\frac{k}{n}\left(1 - \frac{k}{n}\right)\right]^{1-\gamma} \left(\frac{1}{k}\sum_{i=1}^{k} \delta_{X_i} - \frac{1}{n-k}\sum_{i=k+1}^{n} \delta_{X_i}\right),$$

and $\gamma$ is a parameter satisfying $0 \leq \gamma < 1$. The estimator proposed in [6] corresponds to the case of $\gamma = 1/2$.

Estimators of this type consider all possible locations of the change-point, $k$. For each possible $k$ they compute the difference between the empirical probability distributions for the data points on either side of the proposed change-point. This difference is then multiplied by the weighting factor $[k/n(1-k/n)]^{1-\gamma}$. We then require a seminorm, $N$, to measure the difference between the empirical probability distributions. The estimator $\hat{\theta}_n$ is chosen to maximize the difference between the empirical probability distributions under the given seminorm. The weighting factor is required, otherwise values of $k$ near the end points give rise to empirical distributions that contain few data points and therefore give very large statistical errors.

In the theorems stated below we will develop a framework that can deal with a very general class of estimators. Different seminorms represent using different measures of the difference between the distributions before and after the change-point. In the following, we give some examples of seminorms that have been used to estimate change-points for independent data. We will show that these estimators, and a much wider class, are also appropriate for estimating change-points in dependent data. For a measure $\nu$ on $E$ and $f:E \to \mathbb{R}$, we define $\nu(f)$ as

$$\nu(f) \equiv \int f(x)\nu(dx). \tag{2.3}$$

For each choice of seminorm, we require a family of functions that we denote by $\mathcal{F}$. For example, for parametric estimators that only consider a single moment, $\mathcal{F}$ will only contain a single function.

EXAMPLE 1. For a family of functions $\mathcal{F} = \{\mathbb{1}_{.<X_i}, i = 1, \ldots, n\}$, we define norms of a measure $\nu$ via the quantities $d_i = \nu(\mathbb{1}_{.<X_i})$. This corresponds to the setting of [4]. For example,

$$N(\nu) = \sup_{1 \le i \le n} |d_i| \tag{2.4}$$

corresponds to the $L^\infty$ or Kolmogorov–Smirnov norm and

$$N_p(\nu) = \left(\frac{1}{n}\sum_{i=1}^{n}|d_i|^p\right)^{1/p} \tag{2.5}$$

corresponds to the $L^p$ norm. The cases $p=1$ and $p=2$, correspond to the most commonly used $L^1$ and $L^2$ norms. Observe that in this example the family is random and therefore the seminorm is also random.

EXAMPLE 2. For $\mathcal{F} = \{f^p : x \to x^p, p = 1, \ldots, +\infty\}$ we define the seminorm by

$$N(\nu) \equiv \sum_{f \in \mathcal{F}} d(f)|\nu(f)|,$$

where $d(f)$ is a sequence of positive weights. This includes the parametric estimators in which we estimate a change in some moments. For example, differences in the $p$th moment can be detected using the seminorm that applies the measure (2.2) to the function $f^p : x \to x^p$. This framework can also deal with a weighted sum of all moments. This family requires high moments of the marginal law to be finite. To overcome this restriction, one can consider truncated moments, that is, a family given by $\mathcal{F} = \{f_M^p : x \to x^p \mathbb{1}_{|x|<M}, p = 1, \ldots, +\infty\}$, where $M$ is a constant which can be arbitrarily large.

EXAMPLE 3. $\mathcal{F} = \{\mathbb{1}_D, D \in \mathcal{D}\}$, where $\mathcal{D}$ is a family of sets which satisfies certain conditions, such as the family being a VC subclass (see [6]). This means that the family of sets has a covering number which grows polynomially (see [14]).

We now turn our attention to the dependence structure of the sequence. We note that for any given norm, one must apply the measure (2.2) to a family of functions. In this paper we will consider a very general class of dependence structures. For a given sequence we will allow the estimator to use families of functions that satisfy the following condition.

ASSUMPTION 1. There exist constants $C > 0$ and $\rho > 0$ that are independent of the sequence length such that

(2.6) $$\sup_{f \in \mathcal{F}} \sup_{1 \leq i \leq n-m} |\mathrm{corr}(f(X_i), f(X_{i+m}))| \leq Cm^{-\rho}.$$

This assumption simply states that for each of the functions $f$ in $\mathcal{F}$ the correlation between $f(X_i)$ and $f(X_{i+m})$ must decay algebraically or faster with $m$ as $m \to \infty$. This assumption is satisfied for a very general class of data. We now give some examples for which Assumption 1 is satisfied.

EXAMPLE 4. Let $G_1$ and $G_2$ be any measurable functions and $(Z_i)$ be a (possibly nonstationary) Gaussian sequence such that $\sup_{1 \leq i \leq n-m} |\mathrm{corr}(Z_i, Z_{i+m})| \leq Cm^{-\rho}$ and $X_i = G_1(Z_i)$ if $i \leq n\theta$ and $X_i = G_2(Z_i)$ if $i > n\theta$. Then for any family $\mathcal{F}$ such that $\mathbb{E}(f^2(X_i)) < \infty$ for $f \in \mathcal{F}$,

(2.7) $$\sup_{f \in \mathcal{F}} \sup_{1 \leq i \leq n-m} |\mathrm{corr}(f(X_i), f(X_{i+m}))| \leq Cm^{-\rho}$$

(see, e.g., [1]). In fact, this example can be extended to functions of Gaussian vectors using the results of [1].

EXAMPLE 5. Let $(X_i)$ be defined by $X_i \equiv Z_i^{(1)} \equiv \sum_{k=-\infty}^{+\infty} b_k^{(1)} \epsilon_{i-k}^{(1)}$ if $i \leq n\theta$ and $X_i \equiv Z_i^{(2)} \equiv \sum_{k=-\infty}^{+\infty} b_k^{(2)} \epsilon_{i-k}^{(2)}$ if $i > n\theta$, where $(b_k^{(1)})$ and $(b_k^{(2)})$ are real sequences and $(\epsilon_k^{(1)})$ and $(\epsilon_k^{(2)})$ are random stationary sequences with zero mean and finite variance. If $\sum_{k,l=-\infty}^{+\infty} |b_k^{(i)} b_l^{(j)} \mathbb{E}(\epsilon_0^{(i)} \epsilon_{k-l}^{(j)})| < \infty$ for $i, j = 1, 2$, then $(Z_i^{(j)})$ exists almost surely and $\mathbb{E}((Z_i^{(j)})^2) < \infty$. Let $r(k) = \sup_{i,j=1,2} |\mathbb{E}(\epsilon_0^{(i)} \epsilon_k^{(j)})|$. If we assume that $\sup_k |\frac{r(k+m)}{r(k)}| \leq Cm^{-\alpha}$, for $\alpha > 0$, then $|\mathrm{cov}(X_i, X_{i+m})| \leq C'm^{-\alpha}$. This example includes FARIMA processes

with correlated innovations such as GARCH processes. It allows us to model long-range dependence and time-dependent conditional variance. These two features are frequently encountered in financial time series. So, Assumption 1 is satisfied when $\mathcal{F}$ is the set of the identity function.

In Theorem 1 we develop conditions that can deal with countable families of functions and norms that are bounded by weighted moments. In Theorem 2 we consider the case of uncountable families. In this case we need to control the size of the family. This will be done by using covering numbers defined in Assumption 2.

We begin by considering the case where the class of functions $\mathcal{F}$ is countable and the difference between the distributions before and after the change-point may tend to zero as the sequence length, $n$, tends to infinity. This theorem essentially handles the case in which the norm is bounded by a sum of weighted moments and hence includes most commonly used parametric estimators.

For $f$ in $\mathcal{F}$ we set

$$(2.8) \quad \|f\| \equiv \sup_{n \in \mathbb{N}} (P_n(f^2) + Q_n(f^2))^{1/2} = \sup_{n \in \mathbb{N}} (\mathbb{E}_{P_n}[f^2] + \mathbb{E}_{Q_n}[f^2])^{1/2}.$$

THEOREM 1. *Assume that the norm $N$ satisfies*

$$(2.9) \quad N(\nu) \leq \sum_{f \in \mathcal{F}} d(f)|\nu(f)|,$$

*where $\mathcal{F}$ is a countable family of functions satisfying* (2.6) *and $d(f)$ are positive constants such that $\sum_{f \in \mathcal{F}} d(f)\|f\| < \infty$. We assume that there exists a positive sequence $b_n$ such that*

$$(2.10) \quad \mathbb{P}[N(P_n - Q_n) > b_n] \to 1 \qquad as\ n \to \infty.$$

*Let $\bar{\rho} = \min(1 - \epsilon, \rho)$ for any $\epsilon > 0$, where $\rho$ is given in* (2.6). *If*

$$(2.11) \quad b_n^{-1}[n^{-\bar{\rho}/2}(1 + \ln(n)\mathbb{1}_{\gamma - 1 + \bar{\rho}/2 = 0}) + n^{\gamma - 1}] \to 0 \qquad as\ n \to \infty,$$

*then we have*

$$(2.12) \quad \hat{\theta}_n - \theta = O_p(n^{-1} b_n^{-2/\bar{\rho}}).$$

We note that the largest possible value of $\bar{\rho}$ is strictly less than unity and so as long as $\gamma < 1/2$ we will always have $\gamma - 1 + \bar{\rho}/2 \neq 0$, in which case we obtain a less restrictive condition than [6] on the speed at which the difference between the distributions before and after the jump tends to zero. Moreover, if the sequence is LRD ($\rho < 1$), then we have more freedom in the choice of $\gamma$, namely $\gamma \leq 1 - \rho/2$.

This theorem takes a simpler form when $N(P_n - Q_n)$ is bounded away from zero. This is stated in following corollary.

COROLLARY 1. *Under Assumption 1, assume that the seminorm $N$ satisfies* (2.9) *and* (2.10) *with $b_n \geq b > 0$. Then*

$$\hat{\theta}_n - \theta = O_p(n^{-1}). \tag{2.13}$$

Corollary 1 includes the commonly encountered case in which the distributions $P_n$ and $Q_n$ do not depend on the sequence length and the seminorm is nonrandom.

Equation (2.10) controls the rate at which the seminorm of the difference between the two distributions decays to zero by stating that it decays more slowly than some sequence $b_n$. In particular, if the seminorm is nonrandom, one can take $b_n = 2^{-1} N(P_n - Q_n)$. Equation (2.11) requires that random fluctuations arising from sums of the type (2.2), which have size $O(n^{-\rho/2} + n^{\gamma - 1})$, decay to zero faster than the sequence $b_n$ and consequently decay faster than the distance between the two distributions. This is a natural condition to be able to detect a change-point.

We now turn our attention to the case when the family $\mathcal{F}$ contains an uncountable infinity of functions. The following theorem deals with an extremely general set of norms including all of those considered by Carlstein [4]. In this case, under the assumptions that the family has a finite covering number, we obtain the same rate of convergence as in (2.13) when $P_n$ and $Q_n$ are independent of $n$. For the case in which the size of the difference between $P_n$ and $Q_n$ tends to zero as $n \to \infty$ we obtain a rate that depends on the covering number that will typically represent some loss on (2.12).

ASSUMPTION 2. Given two functions $l$ and $u$, the bracket $[l, u]$ is the set of all functions $f$ with $l \leq f \leq u$. Given a norm $\|\cdot\|$ on a space containing $\mathcal{F}$, an $\varepsilon$-bracket for $\|\cdot\|$ is a bracket $[l, u]$ with $\|l - u\| < \varepsilon$. The bracketing number $N_{[\,]}(\varepsilon, \mathcal{F}, \|\cdot\|)$ is the minimal number of $\varepsilon$-brackets needed to cover $\mathcal{F}$.

A family $\mathcal{F}$ is said to satisfy Assumption 2 if

$$\forall \varepsilon > 0 \qquad N_{[\,]}(\varepsilon, \mathcal{F}, \|\cdot\|_X) < \infty, \tag{2.14}$$

where $\|\cdot\|_X$ is a norm satisfying $\sup_{n \in \mathbb{N}} |P_n(|f|)| + |Q_n(|f|)| \leq \|f\|_X$.

We refer the reader to the monograph of van der Vaart and Wellner [14] for examples about bracketing numbers.

The following theorem considers the case when the difference between the distributions before and after the change-point may tend to zero.

THEOREM 2. *Assume that the seminorm satisfies*

$$N(\nu) \leq \sup\{|\nu(f)|, f \in \mathcal{F}\}, \tag{2.15}$$

where $\mathcal{F}$ is a family of functions that satisfies $\sup\{\|f\|, f \in \mathcal{F}\} < \infty$ and Assumptions 1 and 2. Let $\bar{\rho} = \min(1 - \epsilon, \rho)$ for any $\epsilon > 0$, where $\rho$ is given in (2.6), and $\varepsilon_n$ be any positive sequence that tends to zero as $n \to \infty$. We assume that there exists a positive sequence $b_n$ such that

$$\mathbb{P}(N(P_n - Q_n) > b_n) \to 1 \qquad as \ n \to \infty \tag{2.16}$$

and

$$b_n^{-1} N_{[\cdot]}(b_n \varepsilon_n, \mathcal{F}, \|\cdot\|_X)[n^{-\bar{\rho}/2}(1 + \ln(n)\mathbb{1}_{\gamma-1+\bar{\rho}/2=0}) + n^{\gamma-1}] \to 0.$$

Then we have

$$\hat{\theta}_n - \theta = O_p(n^{-1}[b_n^{-1} N_{[\cdot]}(b_n \varepsilon_n, \mathcal{F}, \|\cdot\|_X)]^{2/\bar{\rho}}). \tag{2.17}$$

The following corollary considers the case in which the norm between the distributions before and after the change-point is strictly positive. Provided that the bracketing number is finite, the $n^{-1}$ convergence rate is achieved for any norm within a class of functions satisfying Assumptions 1 and 2.

COROLLARY 2. *Under Assumptions 1 and 2, assume that the seminorm satisfies* (2.10) *with $b_n \geq b > 0$ and* (2.15). *Then* (2.13) *is satisfied.*

REMARK 1. In the case $b_n > b > 0$, Theorems 1 and 2 both give the same $O_p(n^{-1})$ rate for both $\rho < 1$ and $\rho \geq 1$. For Theorem 1, in the case $b_n \to 0$ with $\rho \geq 1$, it is possible to obtain the rate $O_p(n^{-1} b_n^{-2} \ln^2(nb_n^2))$ which can represent a marginally better result. A similar result can be obtained for Theorem 2 with $b_n \to 0$ and $\rho \geq 1$. These results can be obtained by modifying Lemma 1 of our proof using Theorem 3 in [13].

REMARK 2. Assumption 1 can be replaced by the following more general, but less intuitive, condition: there exist constants $C > 0$ and $\rho > 0$, such that for any $m$

$$\sup_{f \in \mathcal{F}} \sup_{k,m,k+m \leq n} \mathbb{E}\left(\sum_{i=k}^{k+m} [f(X_i) - \mathbb{E}(f(X_i))]\right)^2 \leq C m^{2-\rho}. \tag{2.18}$$

In this case $\|f\|$ can be replaced by unity in the assertions of Theorems 1 and 2. Observe that this assumption is particularly weak and satisfied by a large class of processes and families of functions. We now present more examples of commonly used time series models and families of functions that satisfy (2.18) and Assumption 2.

EXAMPLE 6. We begin by considering a linear process with a family of functions that satisfies a Lipschitz condition. Let $\mathcal{F}$ be a family of uniformly bounded functions such that $\sup_{f \in \mathcal{F}} |f(x) - f(y)| \leq C_1 |x-y|^{\eta_1}$ for some $\eta_1 > 0$ and $C_1 > 0$. Then according to [14] $\mathcal{F}$ satisfies Assumption 2 for any $L^p$ norm. We now show that if the sequence is drawn from Example 5, then (2.18) is satisfied under additional weak conditions. Let $X_i^v \equiv \sum_{|k|<v} b_k^{(j)} \epsilon_{i-k}^{(j)}$ with $j=1$ if $i \leq n\theta$ and $j=2$ if $i > n\theta$. Assume that $(\epsilon_k^{(1)}, \epsilon_k^{(2)})$ are $q$-dependent and

(2.19) $$\exists \eta_2 > 0 \quad \forall v \quad \mathbb{E}[X_i - X_i^v]^2 \leq C_2 v^{-\eta_2}.$$

For example, if $|b_k^{(1)}| + |b_k^{(2)}| \leq C|k|^{-\beta}$ and $\beta > 1/2$, then one can readily show that (2.19) is satisfied. The sequence $X_i^v$ is $3v$-dependent for $v > q$, and so by using a blocking technique we have $\mathbb{E}(\sum_{i=k}^{k+m} \bar{f}(X_i^v))^2 \leq Cmv$, where $\bar{f}(X) = f(X) - \mathbb{E}[f(X)]$. Letting $v = m^{1/(1+\eta_1 \eta_2)}$, we obtain

$$\mathbb{E}\left(\sum_{i=k}^{k+m} \bar{f}(X_i)\right)^2 \leq 2\mathbb{E}\left(\sum_{i=k}^{k+m} \bar{f}(X_i^{v(m)})\right)^2 + 2\mathbb{E}\left(\sum_{i=k}^{k+m} (\bar{f}(X_i) - \bar{f}(X_i^{v(m)}))\right)^2$$
$$\leq Cm^{2-\eta_1 \eta_2/(1+\eta_1 \eta_2)}.$$

So (2.18) is also satisfied and hence Theorem 2 applies.

EXAMPLE 7. In this example we consider a linear process given in Example 5 with a family composed of indicator functions, namely $\mathcal{F} = \{f_x(\cdot) \equiv \mathbf{1}_{\cdot \leq x}, x \in \mathbb{R}\}$. This family is relevant to the commonly used $L^p$ and $L^\infty$ norms in Example 1 for which Assumption 2 is satisfied. We assume that $(\epsilon_k^{(1)}, \epsilon_k^{(2)})$ are $q$-dependent and $|b_k^{(1)}| + |b_k^{(2)}| \leq C|k|^{-\beta}$ with $\beta > 1/2$. We begin by assuming that $q = 1$. Then we have

$$\mathbb{E}\left(\sum_{i=k}^{k+m} \bar{f}_x(X_i)\right)^2 \leq 2\mathbb{E}\left(\sum_{i=k}^{k+m} \bar{f}_x(X_i^v)\right)^2 + 2m^2 \sup_i \mathbb{E}[\bar{f}_x(X_i) - \bar{f}_x(X_i^v)]^2,$$

where $\bar{f}_x(X) = f_x(X) - \mathbb{E}[f_x(X)]$. Again, using the blocking technique, we have $\mathbb{E}(\sum_{i=k}^{k+m} \bar{f}_x(X_i^v))^2 \leq Cmv$. One can also show that for some $\eta_1 > 0$, $\sup_x \sup_i \mathbb{E}[f_x(X_i) - f_x(X_i^v)]^2 \leq Cv^{-\eta_1}$. Then by choosing $v \sim m^{1/(1+\eta_1)}$ we obtain $\mathbb{E}(\sum_{i=k}^{k+m} \bar{f}_x(X_i))^2 \leq Cm^{2-\eta_1/(1+\eta_1)}$. The case of $q > 1$ can be handled by dividing the sum $\sum_{i=k}^{k+m} \bar{f}_x(X_i)$ into $q$ blocks such that within each block the innovations are independent. Hence (2.18) is satisfied and Theorem 2 applies.

Before presenting the proofs, we give an intuitive explanation of why the rate of convergence of the estimator does not depend on the dependence structure of the sequence. We define $t_k \equiv k/n$. Then $D_k \equiv D_n(t_k)$, where

$$D_n(t) = t^{1-\gamma}(1-t)^{1-\gamma}\left(\frac{1}{nt}\sum_{i=1}^{[nt]} \delta_{X_i} - \frac{1}{n(1-t)}\sum_{i=[nt]+1}^{n} \delta_{X_i}\right)$$

and $w(t) = t^\gamma (1-t)^\gamma$. We rewrite $D_n(t)$ as the sum of its mean and a centered random component, $B_n(t)$,

$$(2.20) \qquad D_n(t) = \frac{1}{w(t)}[(P_n - Q_n)g(t) + B_n(t)],$$

where $g(t) = t(1 - \theta_n)\mathbb{1}_{t \leq \theta_n} + \theta_n(1-t)\mathbb{1}_{t > \theta_n}$ is a piecewise linear function that takes its maximum at the point $\theta_n \equiv [n\theta]/n$ and $B_n$ is the empirical bridge measure given by

$$(2.21) \qquad B_n(t) = W_n(t) - tW_n(1),$$

$$(2.22) \qquad W_n(t) = \frac{1}{n}\sum_{i=1}^{[nt]}[\delta_{X_i} - \mathscr{L}(X_i)].$$

Our main results stated in Theorems 1 and 2 occur because of the cancellation of two competing effects. One of the effects is concerned with the absolute magnitude of the random noise in $D_n(t)$. The mean component of $D_n(t)$ is monotonically increasing for $t < \theta_n$ and monotonically decreasing for $t > \theta_n$ and therefore takes its maximum at $t = \theta_n$. The estimator is chosen by maximizing $N(D_n(t))$, so if the noise is sufficiently small we would expect to obtain a good estimate. For independent or SRD sequences the partial sums in the centered random component of $D_n(t)$, namely $B_n(t)$, typically have a magnitude of order $n^{-1/2}$ as $n \to \infty$. As shown by Dümbgen, this gives rise to typical errors of order $n^{-1}$ in the estimator. For LRD sequences the partial sums decay more slowly. This means that the stronger the dependence the larger the random component in (2.2). This effect makes the estimation more difficult. One might naively expect that this would mean that LRD sequences have a slower rate of convergence than SRD or independent sequences. However, there is another effect that is concerned with the variations in the noise in the vicinity of the change-point. Correlations in LRD sequences imply that the random noise $B_n(t)$ becomes correlated. This means that the random noise has less rapid variation and local fluctuations become smaller. Estimation requires one to find the global maximum of $N(D_n(t))$ and this depends critically on the local variations in the vicinity of the change-point rather than on the absolute magnitude of the noise. Hence the smaller the local fluctuations are, the easier the estimation becomes. These two effects exactly compensate and give the surprising feature that the overall rate of convergence is the same for all dependence structures.

**3. Simulations.** In this section we present the results of numerical simulations that investigate some of the important practical features of change-point estimation. We confirm that the rate of convergence is $O_p(n^{-1})$ for LRD, SRD and independent sequences. We also determine how large the sequence length needs to be before the $O_p(n^{-1})$ rate is observed.

We considered the estimation of the change-point for a sequence that is a function of a dependent Gaussian variable, $(Y_i)_{i=1,\ldots,n}$ with zero mean and unit variance. We generated a sequence with a change in the marginal distribution by taking

$$X_i = \begin{cases} Y_i^2 - 1, & \text{if } i \leq n\theta, \\ 1 - Y_i^2, & \text{if } i > n\theta. \end{cases}$$

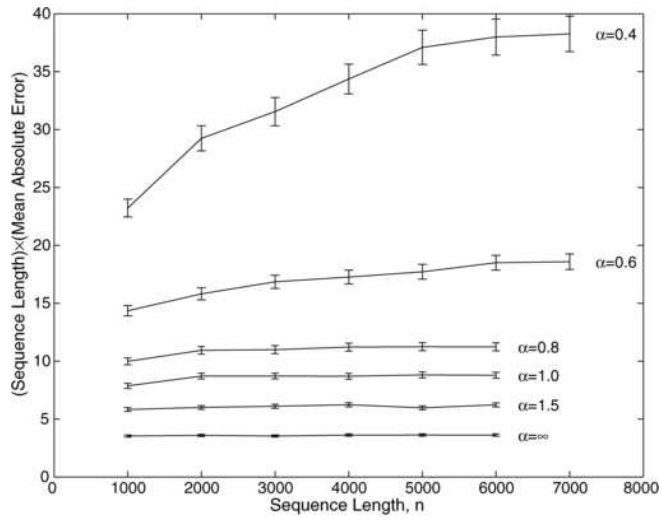

(a)

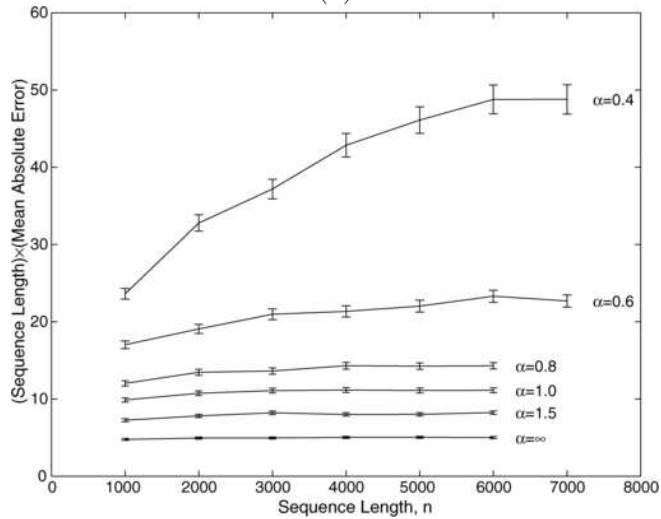

(b)

FIG. 1. *The MAE of $n(\hat{\theta}_n - \theta)$ for different values of $\alpha$: under* (a) *the $L^1$ norm and* (b) *the KS norm.*

The sequence $(X_i)$ has the property that the marginal distributions before and after the jump have the same mean and variance, but have different skewness. We generated the Gaussian sequences $(Y_i)$ with a covariance given by $r(n) = (1+n^2)^{(-\alpha/4)} \sim n^{-\alpha/2}$ using the Durbin–Levinson algorithm (see, e.g., [3]). The sequence $(X_i)$ satisfies Assumption 1 with $\rho = \alpha$. We note that the Durbin–Levinson algorithm has complexity $O(n^2)$ and so generating long sequences can be quite computationally expensive.

We show results for the estimator that uses the Kolmogorov–Smirnov norm (KS) (2.4) and the $L^1$ norm defined in (2.5) with $p = 1$. The parameter $\gamma$ is equal to 0.5. We note, however, that taking different norms, such as $p = 2$ in equation (2.5), yields qualitatively similar results. We considered independent sequences, SRD sequences with $\alpha = 1.5$ and LRD sequences with values of $\alpha = 1.0, 0.8, 0.6$ and $0.4$. We present simulations in which the sequence length, $n$, varies between 1000 and 7000. The mean absolute error $MAE \equiv \mathbb{E}(|(\hat{\theta}_n - \theta)|)$ for each value of $\alpha$ was estimated using 10,000 different sequences. In Figure 1 we plot $n(MAE)$ against $n$ with 95% confidence intervals. Since $\hat{\theta}_n - \theta = O_p(n^{-1})$, we anticipate that $n(MAE)$ should tend to a constant as $n$ tends to infinity. This is clearly seen in Figure 1 for independent, SRD and LRD sequences. As the range of dependence becomes longer, the value of $n$ required to obtain the $O_p(n^{-1})$ scaling becomes larger. This is because the leading order correction to the $O_p(n^{-1})$ rate contains partial sums that are a factor $n^{-\alpha/2}$ smaller than the leading term. So for small $\alpha$, large values of $n$ are required for the leading order term to dominate the corrections.

**4. Proofs.** We will begin by proving that the estimators are consistent. For Theorem 1 this is straightforward, but for Theorem 2 we require a projection argument to deal with the uncountable size of the family $\mathcal{F}$. Having proved consistency, we then turn our attention to the rates proofs. The rates proofs follow a similar pattern to the consistency proof and the techniques used are similar. In the proofs, $C, C_1, C_2, \ldots$ denote generic constants that are independent of $n$ for $n$ large enough whose values may differ in different equations. In general, $\theta \notin \{k/n : k = 1, \ldots, n\}$ so we have defined $\theta_n \equiv [n\theta]/n$. To prove Proposition 1 below and Theorems 1 and 2 it suffices to prove the assertions with $\theta$ replaced by $\theta_n$. In all of the proofs, we will assume $\rho < 1$ since the proofs can be easily adapted for the case $\rho \geq 1$ by replacing $\rho$ with $\bar{\rho}$.

We require the following lemmas for the proofs. The first one is a maximal inequality which is a special case of Theorem 1 in [13].

LEMMA 1. *Assume* (2.6) *with* $\rho < 1$. *Then there exists a constant* $D(\rho) > 0$ *such that*

$$\mathbb{E}\left(\max_{1\leq k\leq n}\left|\sum_{i=1}^{k}[f(X_i) - \mathbb{E}(f(X_i))]\right|\right)^2 \leq D^2(\rho)\|f\|^2 n^{2-\rho}. \quad (4.1)$$

The second lemma controls the size of the empirical bridge and is a simple consequence of (4.1).

LEMMA 2. *Assume* (2.6) *holds with* $0 < \rho < 1$. *Then there exists a constant* $D(\rho) > 0$ *such that for any* $0 < \kappa \leq 1$

$$\mathbb{E}\left[\sup_{|t-\theta_n|\leq\kappa}|(W_n(t) - W_n(\theta_n))(f)|\right] \leq D(\rho)\|f\|n^{-\rho/2}\kappa^{1-\rho/2} \quad (4.2)$$

*and*

$$\mathbb{E}\left[\sup_{|t|\leq\kappa}|(W_n(t))(f)| + \sup_{|t|\leq\kappa}|(B_n(t))(f)|\right] \leq D(\rho)\|f\|(n^{-\rho/2}\kappa^{1-\rho/2}). \quad (4.3)$$

The third lemma controls the size of oscillations of the weighted empirical bridge which we define as

$$B_n^w(t) \equiv w^{-1}(t)B_n(t).$$

LEMMA 3. *Assume* (2.6) *with* $\rho < 1$. *Then there exist constants* $C(\theta, \eta)$ *and* $D(\rho)$ *such that for* $\kappa < \eta$,

$$\mathbb{E}\left(\sup_{|t-\theta_n|\leq\kappa}|(B_n^w(t) - B_n^w(\theta_n))(f)|\right) \leq C(\theta,\eta)D(\rho)\|f\|n^{-\rho/2}\kappa^{1-\rho/2}. \quad (4.4)$$

PROOF. Using Taylor's theorem to expand $w^{-1}(t)$ near $t = \theta_n$, we obtain

$$B_n^w(t) - B_n^w(\theta_n)$$
$$= w^{-1}(\theta_n)(W_n(t) - W_n(\theta_n))$$
$$\quad - (t - \theta_n)[w^{-1}(\theta_n)W_n(1) + (w^{-2}(\xi)w'(\xi))(W_n(t) - tW_n(1))],$$

where $\xi \in (t, \theta_n)$. Therefore, for $\eta$ small enough and $|t - \theta_n| \leq \eta$, there exists a constant $C(\theta, \eta)$ such that

$$\begin{aligned}(4.5)\quad &|(B_n^w(t) - B_n^w(\theta_n))(f)| \\ &\leq w^{-1}(\theta_n)|(W_n(t) - W_n(\theta_n))(f)| + C(\theta,\eta)|t-\theta_n|\sup_{0\leq t\leq 1}|W_n(t)(f)|.\end{aligned}$$

Hence it suffices to control the size of the oscillations of $W_n(t)$. By (4.2) and (4.3) of Lemma 2, we have

$$\mathbb{E}\left(\sup_{|t-\theta_n|\leq\kappa}|(B_n^w(t) - B_n^w(\theta_n))(f)|\right)$$

$$\leq w^{-1}(\theta_n)\mathbb{E}\left(\sup_{|t-\theta_n|\leq\kappa}|(W_n(t) - W_n(\theta_n))(f)|\right)$$

$$+ C(\theta,\eta)\kappa\mathbb{E}\left(\sup_{0\leq t\leq 1}|W_n(t)(f)|\right)$$

$$\leq w^{-1}(\theta_n)D(\rho)\|f\|n^{-\rho/2}\kappa^{1-\rho/2} + D(\rho)\|f\|C(\theta,\eta)\kappa n^{-\rho/2}$$

$$\leq C(\theta,\eta)D(\rho)\|f\|n^{-\rho/2}\kappa^{1-\rho/2},$$

where $C(\theta,\eta)$ may change in each occurrence, and the relation (4.4) follows. □

4.1. *Consistency proofs.* We first recall some notation and introduce some additionally. Let $\delta_n = P_n - Q_n$,

$$h(t) = w^{-1}(t)(t(1-\theta_n)\mathbb{1}_{t\leq\theta_n} + \theta_n(1-t)\mathbb{1}_{t>\theta_n})$$

and $B_n^w(t) = w^{-1}(t)B_n(t)$, where $B_n(t)$ is defined in (2.21) and $w(t) = t^\gamma(1-t)^\gamma$. For $t$ in $G_n \equiv \{k/n, 1\leq k < n\}$ we rewrite $D_n(t)$ defined in (2.20) as

$$D_n(t) = B_n^w(t) + h(t)\delta_n.$$

We also recall that $\hat{\theta}_n$ is a maximum of $\{N(D_n(t)), t \in G_n\}$. The following proposition states the consistency of the estimators.

PROPOSITION 1. *Let $X$ be a sequence and $\mathcal{F}$ a family such that* (2.6) *is satisfied. Assume that the conditions of Theorem* 1 *or Theorem* 2 *are satisfied. Then*

$$\forall \eta > 0 \qquad \mathbb{P}(|\hat{\theta}_n - \theta_n| > \eta) \to 0 \qquad as\ n \to \infty.$$

PROOF OF THEOREM 1. By definition $\hat{\theta}_n$ is a maximum of $N(D_n(t))$. So

(4.6) $$N(D_n(\hat{\theta}_n)) \geq N(D_n(\theta_n)).$$

Using (2.20), we obtain

$$N(B_n^w(\hat{\theta}_n) + \delta_n h(\hat{\theta}_n)) \geq N(B_n^w(\theta_n) + \delta_n h(\theta_n)).$$

Repeated use of the triangle inequality yields

$$N(B_n^w(\hat{\theta}_n)) \geq N(B_n^w(\theta_n) + \delta_n h(\theta_n)) - N(\delta_n h(\hat{\theta}_n))$$

$$\geq N(\delta_n h(\theta_n)) - N(\delta_n h(\hat{\theta}_n)) - N(B_n^w(\theta_n)).$$

Hence,

(4.7) $$N(B_n^w(\hat{\theta}_n)) + N(B_n^w(\theta_n)) \geq N(\delta_n)(h(\theta_n) - h(\hat{\theta}_n)).$$

We define $a_n = \inf_{|t-\theta_n|>\eta}\{h(\theta_n) - h(t)\}$. Then $a_n > a > 0$ for $n$ large enough, because $h$ is monotonically increasing for $t < \theta_n$ and monotonically decreasing for $t > \theta_n$. Since $a_n$ is defined to be an infimum we obtain

(4.8)
$$\begin{aligned}
&\mathbb{P}[|\hat{\theta}_n - \theta_n| > \eta] \\
&= \mathbb{P}[N(B_n^w(\hat{\theta}_n)) + N(B_n^w(\theta_n)) \geq aN(\delta_n), |\hat{\theta}_n - \theta_n| > \eta] \\
&\leq \mathbb{P}[N(B_n^w(\hat{\theta}_n)) + N(B_n^w(\theta_n)) \geq ab_n, |\hat{\theta}_n - \theta_n| > \eta] \\
&\quad + \mathbb{P}[N(\delta_n) \leq b_n].
\end{aligned}$$

We use the fact that $\mathbb{P}[X+Y \geq \varepsilon, B] \leq \mathbb{P}[|X| \geq \varepsilon/2, B] + \mathbb{P}[|Y| \geq \varepsilon/2, B]$, for any random variables $X$ and $Y$, set $B$ and $\varepsilon > 0$, to obtain

(4.9)
$$\begin{aligned}
\mathbb{P}[|\hat{\theta}_n - \theta_n| > \eta] &\leq \mathbb{P}\left[N(B_n^w(\hat{\theta}_n)) \geq \frac{ab_n}{2}, |\hat{\theta}_n - \theta_n| > \eta\right] \\
&\quad + \mathbb{P}\left[N(B_n^w(\theta_n)) \geq \frac{ab_n}{2}\right] + \mathbb{P}[N(\delta_n) \leq b_n] \\
&\equiv A_1 + A_2 + A_3.
\end{aligned}$$

We begin by controlling $A_1$. We will assume that $\eta$ is sufficiently small such that $\theta_n - \eta > 0$ and $1 - \theta_n - \eta > 0$, since other cases can be dealt with similarly. For the sake of brevity we introduce the notation $\beta_{min} = \min(\theta_n - \eta, 1 - \theta_n - \eta)$ and $\beta_{max} = \max(\theta_n - \eta, 1 - \theta_n - \eta)$. We introduce sets $S_1, \ldots, S_J$ given by

$$S_j = \{t : 2^{-j} \leq t(\theta_n - \eta)^{-1} < 2^{-j+1}\} \cup \{t : 2^{-j} \leq (1-t)(1-\theta_n-\eta)^{-1} < 2^{-j+1}\}.$$

The integer $J$ is chosen so that $n^{-1}2^{J-1} \leq \beta_{\max} < n^{-1}2^J$. As $j$ increases these sets become increasingly close to the end points of the domain and $J$ is chosen to be large enough so that the smallest and largest possible values of the change point (i.e., $\hat{\theta}_n = 1/n$ and $\hat{\theta}_n = 1 - 1/n$) are included in one of the sets. Then

(4.10) $$A_1 = \sum_{j=1}^{J} \mathbb{P}\left[\hat{\theta}_n \in S_j, N(B_n^w(\hat{\theta}_n)) \geq \frac{ab_n}{2}\right] \equiv \sum_{j=1}^{J} \tilde{A}_1(n,j),$$

where

(4.11) $$\tilde{A}_1(n,j) \leq \mathbb{P}\left[\sup_{t \in S_j} N(B_n(t)) \geq \frac{ab_n}{2} \inf_{t \in S_j} w(t)\right].$$

A simple calculation shows that $\inf_{t \in S_j} w(t) = \min(w((\theta_n - \eta)2^{-j}),$ $w((1 - \theta_n - \eta)2^{-j}) > \beta_{\min}^{\gamma} 2^{-j\gamma - 1}$. Hence applying the Markov inequality to (4.11) we obtain

$$\tilde{A}_1(n,j) \leq \beta_{\min}^{-\gamma} 2^{j\gamma + 2} a^{-1} b_n^{-1} \mathbb{E}\left[\sup_{t \in S_j} N(B_n(t))\right].$$

In order to control $\mathbb{E}[\sup_{t \in S_j} N(B_n(t))]$ we need to control $\mathbb{E}[\sup_{t \in S_j} |B_n(t)(f)|]$ for $f \in \mathcal{F}$. We use (2.9) to prove Proposition 1 under the conditions of Theorem 1 and use a chaining argument for Proposition 1 under the conditions of Theorem 2. The control of $\mathbb{E}[\sup_{t \in S_j} |B_n(t)(f)|]$ is formulated in Lemma 2. Using (2.9) and applying Lemma 2, we obtain

$$\tilde{A}_1(n,j) \leq \beta_{\min}^{-\gamma} 2^{j\gamma + 2} a^{-1} b_n^{-1} \mathbb{E}\left[\sup_{t \in S_j} \sum_{f \in \mathcal{F}} d(f)|B_n(t)(f)|\right]$$

$$\leq \beta_{\min}^{-\gamma} 2^{j\gamma + 3} a^{-1} b_n^{-1} \sum_{f \in \mathcal{F}} d(f) D(\rho) \|f\| n^{-\rho/2} (\beta_{\max} 2^{-j})^{1-\rho/2}.$$

Substituting the above inequality into (4.10), we obtain

$$(4.12) \quad A_1 \leq 8\beta_{\min}^{-\gamma} \beta_{\max}^{1-\rho/2} a^{-1} b_n^{-1} D(\rho) n^{-\rho/2} \sum_{f \in \mathcal{F}} d(f) \|f\| \sum_{j=1}^{J} 2^{(\gamma - 1 + \rho/2)j}.$$

It is easy to show that

$$(4.13) \quad \sum_{j=1}^{J} 2^{(\gamma - 1 + \rho/2)j}$$
$$\leq C(\rho, \gamma)(1 + n^{\gamma - 1 + \rho/2} \mathbb{1}_{\gamma - 1 + \rho/2 \neq 0} + \ln n \mathbb{1}_{\gamma - 1 + \rho/2 = 0}).$$

Substituting (4.13) into (4.12) and relabeling the constant yields

$$(4.14) \quad A_1 \leq C_1 b_n^{-1}(n^{-\rho/2}(1 + \ln n \mathbb{1}_{\gamma - 1 + \rho/2 = 0}) + n^{\gamma - 1}) \sum_{f \in \mathcal{F}} d(f) \|f\|.$$

To control $A_2$ we make similar use of the Lemma 2 to obtain

$$(4.15) \quad A_2 \leq C_2 D(\rho) b_n^{-1} n^{-\rho/2} \sum_{f \in \mathcal{F}} d(f) \|f\|.$$

Finally from (4.9), (4.14) and (4.15) we deduce

$$\mathbb{P}[|\hat{\theta}_n - \theta_n| > \eta]$$
$$\leq C b_n^{-1} \sum_{f \in \mathcal{F}} d(f) D(\rho) \|f\| (n^{\gamma - 1}(1 + \ln n \mathbb{1}_{\gamma - 1 + \rho/2 = 0}) + n^{-\rho/2})$$
$$+ \mathbb{P}(N(\delta_n) \leq b_n).$$

Taking the limit $n \to \infty$ under the conditions (2.10) and (2.11) and the condition $\sum_{f \in \mathcal{F}} d(f) \|f\| < \infty$ completes the proof. □

PROOF OF THEOREM 2. The consistency proof under the assumptions of Theorem 2 is identical to that of Theorem 1 up until (4.9). Then we proceed by using a projection argument to bound $A_1, A_2$, and $A_3$. This projection argument is to deal with the uncountable family of functions. Since $N(K) \equiv N_{[\cdot]}(2^{-K}, \mathcal{F}, \|\cdot\|_X)$ is finite for any integer $K$, there exists a finite sequence of pairs of functions $(f_i^K, \Delta_i^K)_{1 \leq i \leq N(K)}$, such that $\forall f \in \mathcal{F}$ there exists $i$ such that $|f - f_i^K| \leq \Delta_i^K$, and $\|\Delta_i^K\|_X \leq 2^{-K}$. For each $K$ we define a map $\mathcal{M}$ from $\mathcal{F}$ to $\mathcal{F} \times \mathcal{F}$ by $\mathcal{M}(f) = (f_{i(f)}^K, \Delta_{i(f)}^K) \equiv (\pi_K(f), \Delta_K(f))$, where $i(f) = \inf\{1 \leq i \leq N(K) | f_i^K - \Delta_i^K \leq f \leq f_i^K + \Delta_i^K\}$.

We *assume* that $\gamma - 1 + \rho/2 \neq 0$ (the case $\gamma - 1 + \rho/2 = 0$ can be handled similarly and is hence omitted). We apply the Markov inequality to $A_1$ in equation (4.9) and then use the assumption (2.15) on the seminorm $N$ to obtain

$$A_1 \leq 2a^{-1}b_n^{-1} \mathbb{E}\left(\sup_{f \in \mathcal{F}} |B_n^w(\hat{\theta}_n)(f)| \mathbb{1}_{|\hat{\theta}_n - \theta_n| > \eta}\right).$$

To control $A_1$ we will consider two cases: $\hat{\theta}_n > \theta_n$ and $\hat{\theta}_n < \theta_n$, hence

$$A_1 \leq 2a^{-1}b_n^{-1}\bigg(\mathbb{E}\left(\sup_{f \in \mathcal{F}} |B_n^w(\hat{\theta}_n)(f)| \mathbb{1}_{0 < \hat{\theta}_n < \theta_n - \eta}\right)$$

$$+ \mathbb{E}\left(\sup_{f \in \mathcal{F}} |B_n^w(\hat{\theta}_n)(f)| \mathbb{1}_{\theta_n + \eta < \hat{\theta}_n < 1}\right)\bigg)$$

$$\equiv A_1' + A_1''.$$

We first control $A_1'$. Writing $f = f - \pi_K(f) + \pi_K(f)$ gives

(4.16)
$$A_1' \leq 2a^{-1}b_n^{-1} \mathbb{E}\left(\sup_{f \in \mathcal{F}} |B_n^w(\hat{\theta}_n)(f - \pi_K(f))| \mathbb{1}_{0 < \hat{\theta}_n < \theta_n - \eta}\right)$$
$$+ 2a^{-1}b_n^{-1} \mathbb{E}\left(\sup_{f \in \mathcal{F}} |B_n^w(\hat{\theta}_n)(\pi_K(f))| \mathbb{1}_{0 < \hat{\theta}_n < \theta_n - \eta}\right).$$

Using the definitions of $B_n$ and $W_n$, we observe that if $|\phi| \leq g$, then

(4.17) $\quad |B_n(t)(\phi)| \leq |W_n(t)(g)| + |tW_n(1)(g)| + 4\left|t \sup_i \mathbb{E}(|g(X_i)|)\right|.$

Using the fact that $|f - \pi_K(f)| \leq \Delta_K(f)$, applying (4.17) to the first term in (4.16) and the triangle inequality to the second term in (4.16), we obtain

$$A_1' \leq 2a^{-1}b_n^{-1} \mathbb{E}\left(\sup_{f \in \mathcal{F}} |w^{-1}(\hat{\theta}_n) W_n(\hat{\theta}_n)(\Delta_K(f))| \mathbb{1}_{0 < \hat{\theta}_n < \theta_n - \eta}\right)$$

$$+ 2a^{-1}b_n^{-1}\mathbb{E}\bigg(\sup_{f\in\mathcal{F}}|w^{-1}(\hat\theta_n)\hat\theta_n W_n(1)(\Delta_K(f))|\mathbb{1}_{0<\hat\theta_n<\theta_n-\eta}\bigg)$$

$$+ 8a^{-1}\sup_{0<t<\theta_n-\eta}|w^{-1}(t)t|b_n^{-1}\sup_{f\in\mathcal{F}}\sup_{P\in\{P_n,Q_n\}}\mathbb{E}_P(|\Delta_K(f)|)$$

$$+ 2a^{-1}b_n^{-1}\mathbb{E}\bigg(\sup_{f\in\mathcal{F}}|w^{-1}(\hat\theta_n)W_n(\hat\theta_n)(\pi_K(f))|\mathbb{1}_{0<\hat\theta_n<\theta_n-\eta}\bigg)$$

$$+ 2a^{-1}b_n^{-1}\mathbb{E}\bigg(\sup_{f\in\mathcal{F}}|w^{-1}(\hat\theta_n)\hat\theta_n W_n(1)(\pi_K(f))|\mathbb{1}_{0<\hat\theta_n<\theta_n-\eta}\bigg)$$

$$\equiv A'_{1,1} + A'_{1,2} + A'_{1,3} + A'_{1,4} + A'_{1,5}.$$

Following a similar procedure used in the proof of Proposition 1 under the conditions of Theorem 1, we introduce the sets $S'_1,\ldots,S'_{J'}$ defined as

$$S'_j = \{t : 2^{-j} < t{\beta'}^{-1} \le 2^{-j+1}\}.$$

Without loss of generality we assume $\beta' \equiv \theta_n - \eta > 0$ and choose $J'$ to be the integer such that $n^{-1} \in S'_{J'}$, hence $\beta' 2^{-J'} < n^{-1} \le \beta' 2^{-J'+1}$. The proof proceeds in a similar way to that of Proposition 1 under the conditions of Theorem 1. We control $A'_{1,1}$ using Lemma 2 to obtain

$$(4.18) \qquad A'_{1,1} \le C_1 b_n^{-1} \sum_{i=1}^{N(K)} D(\rho)\|\Delta_i^K\|(n^{\gamma-1} + n^{-\rho/2}).$$

Similar use of Lemma 2 on $A'_{1,2}$ yields

$$(4.19) \qquad A'_{1,2} \le C(\gamma,\theta_n,\eta)b_n^{-1}\sum_{i=1}^{N(K)} D(\rho)\|\Delta_i^K\|n^{-\rho/2}.$$

Similar bounds hold for $A'_{1,4}$ and $A'_{1,5}$. Combining these four bounds with the fact that $\sup_{f\in\mathcal{F}}\sup_{P\in\{P_n,Q_n\}}\mathbb{E}_P(|\Delta_K(f)|) \le 2^{-K}$ we obtain

$$(4.20)\quad A'_1 \le C_1 b_n^{-1}\sum_{i=1}^{N(K)} D(\rho)[\|\Delta_i^K\| + \|f_i^K\|](n^{\gamma-1} + n^{-\rho/2}) + C_1 b_n^{-1} 2^{-K}.$$

A similar bound can be derived for $A''_1$. Hence, we conclude that

$$(4.21)\quad A_1 \le C_1 b_n^{-1}\sum_{i=1}^{N(K)} D(\rho)[\|\Delta_i^K\| + \|f_i^K\|](n^{\gamma-1} + n^{-\rho/2}) + C_1 b_n^{-1} 2^{-K},$$

where $C_1$ is a constant that depends only on $\gamma, \theta$ and $\eta$.

To control $A_2$ we write

$$A_2 \leq 2a^{-1}b_n^{-1}\mathbb{E}\left(\sup_{f\in\mathcal{F}}\left|\frac{B_n(\theta_n)}{w(\theta_n)}(f-\pi_K(f))\right|\right)$$
$$+ 2a^{-1}b_n^{-1}\mathbb{E}\left(\sup_{f\in\mathcal{F}}\left|\frac{B_n(\theta_n)}{w(\theta_n)}(\pi_K(f))\right|\right).$$

Applying (4.17) to the first term on the right-hand side of the above equation gives

$$A_2 \leq 2a^{-1}b_n^{-1}\mathbb{E}\left(\sup_{f\in\mathcal{F}}\left|\frac{W_n(\theta_n)}{w(\theta_n)}(\Delta_K(f))\right|\right)$$
$$+ 2a^{-1}b_n^{-1}\mathbb{E}\left(\sup_{f\in\mathcal{F}}\left|\frac{\theta_n W_n(1)}{w(\theta_n)}(\Delta_K(f))\right|\right)$$
$$+ 8\theta_n w^{-1}(\theta_n)a^{-1}b_n^{-1}2^{-K} + 2a^{-1}b_n^{-1}\mathbb{E}\left(\sup_{f\in\mathcal{F}}|B_n^w(\theta_n)(\pi_K(f))|\right).$$

Hence, again by Lemma 2 we have

$$(4.22) \quad A_2 \leq C_2 b_n^{-1}\sum_{i=1}^{N(K)} D(\rho)[\|\Delta_i^K\| + \|f_i^K\|]n^{-\rho/2} + C_2 b_n^{-1} 2^{-K},$$

where $C_2$ is some constant depending only on $\gamma, \theta$ and $\eta$. Finally, from (4.9), (4.21) and (4.22) we have

$$\mathbb{P}(|\hat{\theta}_n - \theta_n| > \eta) \leq C_3 b_n^{-1} D(\rho) N(K) \sup_{1\leq i\leq N(K)}[\|\Delta_i^K\| + \|f_i^K\|](n^{\gamma-1} + n^{-\rho/2})$$
$$+ C_3 b_n^{-1} 2^{-K} + \mathbb{P}(N(\delta_n) \leq b_n).$$

We choose $K$ such that $2^{-K} \sim b_n \varepsilon_n$, where $\varepsilon_n$ is any positive sequence that tends to zero. Since $b_n$ satisfies $b_n^{-1} N_{[\cdot]}(b_n \varepsilon_n, \mathcal{F}, \|\cdot\|_X)(n^{\gamma-1} + n^{-\rho/2}) \to 0$ and $\mathbb{P}(N(\delta_n) \leq b_n) \to 0$, taking the limit as $n \to \infty$ completes the proof. $\square$

REMARK 3. When $b_n > b > 0$, we choose $K$ to be independent of $n$. We let $n$ tend to infinity and $K$ tend to infinity. This completes the proof without posing any restriction on the rate of $N_{[\cdot]}(\varepsilon, \mathcal{F}, \|\cdot\|_X)$.

4.2. *Proof of Theorem* 1. Let $M$ be a positive integer, $b$ and $c$ be positive real numbers and $r_n$ be a positive sequence. We first show that for $n$ large enough

$$(4.23) \quad \mathbb{P}(r_n|\hat{\theta}_n - \theta_n| > 2^M) \leq E_1 + E_2 + E_3 + \mathbb{P}(|\hat{\theta}_n - \theta_n| > \eta),$$

where

$$E_1 \leq \mathbb{P}[r_n^{-1}2^M < |\hat{\theta}_n - \theta_n| \leq \eta, N(B_n^w(\hat{\theta}_n) - B_n^w(\theta_n)) \geq \tilde{C}|\hat{\theta}_n - \theta_n|],$$

(4.24) $\quad E_2 \equiv \mathbb{P}(N(B_n^w(\theta_n)) > c),$

$$E_3 \equiv \mathbb{P}(N(\delta_n) \leq b),$$

$\tilde{C} \equiv C_h(bh(\theta_n) - 2c)$ and $C_h$ is a constant depending only on $\theta$ and $\gamma$.

Recall that $\delta_n = P_n - Q_n$ and $B_n^w(t) = w^{-1}(t)B_n(t)$, where $B_n(t)$ is defined in (2.21). Then

(4.25) $\quad\quad\quad\quad\quad\quad D_n(t) = B_n^w(t) + h(t)\delta_n.$

For all $t$ we have

$$D_n(t) = B_n^w(t) - B_n^w(\theta_n) + B_n^w(\theta_n)\left(1 - \frac{h(t)}{h(\theta_n)}\right) + \frac{h(t)}{h(\theta_n)}D_n(\theta_n).$$

Applying the seminorm and the triangle inequality to the above expression yields

$$N(D_n(t)) \leq N(B_n^w(t) - B_n^w(\theta_n)) + \left(1 - \frac{h(t)}{h(\theta_n)}\right)N(B_n^w(\theta_n))$$
$$+ \frac{h(t)}{h(\theta_n)}N(D_n(\theta_n)).$$

Therefore

(4.26) 
$$N(D_n(t)) - N(D_n(\theta_n)) \leq N(B_n^w(t) - B_n^w(\theta_n))$$
$$+ \left(\frac{h(t)}{h(\theta_n)} - 1\right)(N(D_n(\theta_n)) - N(B_n^w(\theta_n))).$$

Let $\hat{\theta}_n$ be a maximum of $\{N(D_n(t)), t \in G_n\}$, where $G_n \equiv \{k/n, 1 \leq k < n\}$. Since $\hat{\theta}_n$ is a maximum, we have

$$N(B_n^w(\hat{\theta}_n) - B_n^w(\theta_n)) \geq \left(1 - \frac{h(\hat{\theta}_n)}{h(\theta_n)}\right)[N(D_n(\theta_n)) - N(B^w(\theta_n))].$$

Applying the triangle inequality to (4.25) gives $N(D_n(\theta_n)) \geq N(\delta_n h(\theta_n)) - N(B_n^w(\theta_n))$ and therefore we obtain

$$N(B_n^w(\hat{\theta}_n) - B_n^w(\theta_n)) \geq \left(1 - \frac{h(\hat{\theta}_n)}{h(\theta_n)}\right)[N(\delta_n h(\theta_n)) - 2N(B_n^w(\theta_n))].$$

There exists $C_h$ which depends only on $\theta$ and $\gamma$ such that for all $t \in (0,1)$,

$$\left(1 - \frac{h(t)}{h(\theta_n)}\right) \geq C_h|t - \theta_n|.$$

Therefore we obtain

(4.27) $N(B_n^w(\hat{\theta}_n) - B_n^w(\theta_n)) \geq C_h |\hat{\theta}_n - \theta_n|(N(\delta_n h(\theta_n)) - 2N(B_n^w(\theta_n)))$.

For any positive integer $M$ and any positive constants $b$ and $c$, we have

$$\mathbb{P}(r_n|\hat{\theta}_n - \theta_n| > 2^M) \leq \mathbb{P}(r_n^{-1}2^M < |\hat{\theta}_n - \theta_n| \leq \eta, N(\delta_n) > b, N(B_n^w(\theta_n)) \leq c)$$
$$+ \mathbb{P}(N(B_n^w(\theta_n)) > c) + \mathbb{P}(N(\delta_n) \leq b) + \mathbb{P}(|\hat{\theta}_n - \theta_n| > \eta)$$
$$\equiv E_1 + E_2 + E_3 + \mathbb{P}(|\hat{\theta}_n - \theta_n| > \eta).$$

Now from (4.27) we infer that

$$E_1 \leq \mathbb{P}[r_n^{-1}2^M < |\hat{\theta}_n - \theta_n| \leq \eta,$$
$$N(B_n^w(\hat{\theta}_n) - B_n^w(\theta_n)) \geq C_h(bh(\theta_n) - 2c)|\hat{\theta}_n - \theta_n|].$$

This completes the proof of (4.23).

In order to control $E_1$ we define the shells

(4.28) $$S_{n,j} = \{t : 2^j < r_n|t - \theta_n| \leq 2^{j+1}\},$$

where $r_n$ is a positive sequence to be chosen later. Let $0 < \eta < \min(\theta_n, 1 - \theta_n)/2$ and $J \equiv J(n,\eta)$ be chosen such that $2^J < r_n\eta \leq 2^{J+1}$. From the definitions of the shells $S_{n,j}$ and $J$ we obtain

(4.29) $$E_1 \leq \sum_{j=M}^{J} \mathbb{P}[\hat{\theta}_n \in S_{n,j}, N(B_n^w(\hat{\theta}_n) - B_n^w(\theta_n)) \geq \tilde{C}|\hat{\theta}_n - \theta_n|].$$

Now, for $\hat{\theta}_n \in S_{n,j}$, we have $|\hat{\theta}_n - \theta_n| \geq 2^j r_n^{-1}$. Hence using (2.9) and (4.4), we get

(4.30) $$E_1 \leq \tilde{C}^{-1} C(\theta, \eta) \sum_{j=M}^{J} \sum_{f \in \mathcal{F}} d(f)\|f\|D(\rho)2^{(-1/2)j\rho} r_n^{\rho/2} n^{-\rho/2}.$$

For $E_2$, using (2.9) and Lemma 2, we obtain

(4.31) $$E_2 \leq (cw(\theta_n))^{-1} \sum_{f \in \mathcal{F}} d(f)\|f\|D(\rho) n^{-\rho/2}.$$

Now, from (4.23), (4.30) and (4.31) we obtain

$$\mathbb{P}(r_n|\hat{\theta}_n - \theta_n| > 2^M)$$
$$\leq (bh(\theta_n) - 2c)^{-1} C(\theta,\eta) \sum_{j=M}^{J} \sum_{f \in \mathcal{F}} d(f)\|f\|D(\rho) 2^{-1/2 j\rho} r_n^{1/2\rho} n^{-\rho/2}$$
$$+ (cw(\theta_n))^{-1} D(\rho) \sum_{f \in \mathcal{F}} d(f)\|f\| n^{-\rho/2}$$
$$+ \mathbb{P}(N(\delta_n) \leq b) + \mathbb{P}(|\hat{\theta}_n - \theta_n| > \eta).$$

This inequality holds for any $b, c$ and $r_n$, so choosing $b = b_n$, $c = h(\theta_n)b_n/4$ and $r_n = nb_n^{2/\rho}$ and relabeling the constants yields

$$\mathbb{P}(r_n|\hat{\theta}_n - \theta_n| > 2^M) \leq C(\theta, \eta)D(\rho) \sum_{f \in \mathcal{F}} \sum_{j=M}^{J} d(f)\|f\|2^{-1/2j\rho}$$

$$+ C(\theta)D(\rho) \sum_{f \in \mathcal{F}} d(f)\|f\|b_n^{-1}n^{-\rho/2}$$

$$+ \mathbb{P}(|\hat{\theta}_n - \theta_n| > \eta) + \mathbb{P}(N(\delta_n) \leq b_n).$$

Finally letting $n$, then $M$ tend to infinity completes the proof of Theorem 1.

4.3. *Proof of Theorem* 2. The proof of Theorem 2 is identical to that of Theorem 1 up until (4.23). We proceed by using a projection argument to bound $E_1, E_2$ and $E_3$. From (2.15) and (4.24) we have

$$E_1 \leq \mathbb{P}\left[r_n^{-1}2^M < |\hat{\theta}_n - \theta_n| \leq \eta, \sup_{f \in \mathcal{F}}|(B_n^w(\hat{\theta}_n) - B_n^w(\theta_n))(f)| \geq \tilde{C}|\hat{\theta}_n - \theta_n|\right]$$

$$\leq \tilde{C}^{-1}\mathbb{E}\left[|\hat{\theta}_n - \theta_n|^{-1}\mathbb{1}_{r_n^{-1}2^M < |\hat{\theta}_n - \theta_n| \leq \eta} \sup_{f \in \mathcal{F}}|(B_n^w(\hat{\theta}_n) - B_n^w(\theta_n))(f)|\right].$$

Then from (4.5) we obtain

$$E_1 \leq \tilde{C}^{-1}\mathbb{E}\bigg[|\hat{\theta}_n - \theta_n|^{-1}\mathbb{1}_{r_n^{-1}2^M < |\hat{\theta}_n - \theta_n| \leq \eta}$$

(4.32)
$$\times \sup_{f \in \mathcal{F}}\left|\frac{1}{w(\theta_n)}|(W_n(\hat{\theta}_n) - W_n(\theta_n))(f)|\right|\bigg]$$

$$+ \tilde{C}^{-1}\mathbb{E}\left[\sup_{f \in \mathcal{F}}\left|C(\theta, \eta) \sup_{0 \leq t \leq 1}|W_n(t)(f)|\right|\right] \equiv F_1 + G_1.$$

Using the same projection as in the consistency proof we obtain

$$F_1 \leq \tilde{C}^{-1}\mathbb{E}\bigg[|\hat{\theta}_n - \theta_n|^{-1}\mathbb{1}_{r_n^{-1}2^M < |\hat{\theta}_n - \theta_n| \leq \eta}$$

$$\times \sup_{f \in \mathcal{F}}\left|\frac{1}{w(\theta_n)}(W_n(\hat{\theta}_n) - W_n(\theta_n))(f - \pi_K(f))\right|\bigg]$$

$$+ \tilde{C}^{-1}\mathbb{E}\bigg[|\hat{\theta}_n - \theta_n|^{-1}\mathbb{1}_{r_n^{-1}2^M < |\hat{\theta}_n - \theta_n| \leq \eta}$$

$$\times \sup_{f \in \mathcal{F}}\left|\frac{1}{w(\theta_n)}(W_n(\hat{\theta}_n) - W_n(\theta_n))(\pi_K(f))\right|\bigg].$$

We observe that for any $\phi$ and $g$ such that $|\phi| \leq g$,
$$|(W_n(t) - W_n(\theta_n))(\phi)| \leq |(W_n(t) - W_n(\theta_n))(g)| \\ + 2(|t - \theta_n| + n^{-1}) \sup_i \mathbb{E}(g(X_i)).$$

Since $|f - \pi_K(f)| \leq \Delta_K(f)$, by choosing $\phi = f - \pi_K(f)$ and $g = \Delta_K(f)$ in the above inequality, we have

$$F_1 \leq \tilde{C}^{-1} \mathbb{E}\bigg[|\hat{\theta}_n - \theta_n|^{-1} \mathbb{1}_{r_n^{-1} 2^M < |\hat{\theta}_n - \theta_n| \leq \eta} \\ \times \sup_{f \in \mathcal{F}} \bigg| \frac{1}{w(\theta_n)} (W_n(\hat{\theta}_n) - W_n(\theta_n))(\Delta_K(f)) \bigg| \bigg] \\ + 2\tilde{C}^{-1}(1 + n^{-1} r_n 2^{-M}) \sup_{f \in \mathcal{F}} \sup_{P \in \{P_n, Q_n\}} \mathbb{E}_P(|\Delta_K(f)|) \\ + \tilde{C}^{-1} \bigg[|\hat{\theta}_n - \theta_n|^{-1} \mathbb{1}_{r_n^{-1} 2^M < |\hat{\theta}_n - \theta_n| \leq \eta} \\ \times \sup_{f \in \mathcal{F}} \bigg| \frac{1}{w(\theta_n)} (W_n(\hat{\theta}_n) - W_n(\theta_n))(\pi_K(f)) \bigg| \bigg] \\ \equiv F_{1,1} + F_{1,2} + F_{1,3}.$$

Using the decomposition of $\{t : r_n^{-1} 2^M < |\hat{\theta}_n - \theta_n| \leq \eta\}$ over the shells defined in (4.28), we obtain

$$F_{1,1} \leq \tilde{C}^{-1} \sum_{j=M}^{J} \mathbb{E}\bigg[|\hat{\theta}_n - \theta_n|^{-1} \mathbb{1}_{\hat{\theta}_n \in S_{n,j}} \sup_{f \in \mathcal{F}} \bigg| \frac{1}{w(\theta_n)} (W_n(\hat{\theta}_n) - W_n(\theta_n))(\Delta_K(f)) \bigg| \bigg] \\ \leq \tilde{C}^{-1} \sum_{j=M}^{J} \mathbb{E}\bigg[(2^j r_n^{-1})^{-1} \sup_{f \in \mathcal{F}} \bigg| \frac{1}{w(\theta_n)} (W_n(\hat{\theta}_n) - W_n(\theta_n))(\Delta_K(f)) \bigg| \bigg] \\ \leq \tilde{C}^{-1} \sum_{i=1}^{N(K)} \sum_{j=M}^{J} (2^j r_n^{-1})^{-1} \mathbb{E}\bigg[\sup_{t \in S_{n,j}} \bigg| \frac{1}{w(\theta_n)} (W_n(t) - W_n(\theta_n))(\Delta_i^K) \bigg| \bigg].$$

By (4.2) of Lemma 2 we get

(4.33) $\quad F_{1,1} \leq \tilde{C}^{-1} w^{-1}(\theta_n) D(\rho) \sum_{i=1}^{N(K)} \sum_{j=M}^{J} \|\Delta_i^K\| n^{-\rho/2} 2^{-(1/2)j\rho+1} r_n^{\rho/2}.$

A similar bound holds for $F_{1,3}$, and since $\sup_{f \in \mathcal{F}} \sup_{P \in \{P_n, Q_n\}} \mathbb{E}_P(|\Delta_K(f)|) \leq 2^{-K}$ we get

$$F_1 \leq \tilde{C}^{-1} \left(\frac{r_n}{n}\right)^{\rho/2} w^{-1}(\theta_n) D(\rho) \sum_{i=1}^{N(K)} \sum_{j=M}^{J} [\|\Delta_i^K\| + \|f_i^K\|] 2^{(-1/2)j\rho+1}$$

(4.34)
$$+ \tilde{C}^{-1}(1 + n^{-1}r_n 2^{-M})2^{-K}.$$

Similarly, we have

(4.35)
$$G_1 \leq (bh(\theta_n) - 2c)^{-1} C(\theta, \eta) D(\rho) \sum_{i=1}^{N(K)} [\|\Delta_i^K\| + \|f_i^K\|] n^{-\rho/2}$$
$$+ (bh(\theta_n) - 2c)^{-1} C(\theta, \eta) 2^{-K}.$$

For $r_n \leq n$, we have $n^{-1} r_n 2^{-M} \leq 1$. From (4.32), (4.34) and (4.35) we obtain

(4.36)
$$E_1 \leq (C_h(bh(\theta_n) - 2c))^{-1} \left(\frac{r_n}{n}\right)^{\rho/2} w^{-1}(\theta_n) D(\rho)$$
$$\times \sum_{i=1}^{N(K)} \sum_{j=M}^{J} [\|\Delta_i^K\| + \|f_i^K\|] 2^{(-1/2)j\rho + 1}$$
$$+ (bh(\theta_n) - 2c)^{-1} C(\theta, \eta) D(\rho) \sum_{i=1}^{N(K)} [\|\Delta_i^K\| + \|f_i^K\|] n^{-\rho/2}$$
$$+ (bh(\theta_n) - 2c)^{-1} C(\theta, \eta) 2^{-K}.$$

For $E_2$ we use a similar argument to obtain

(4.37) $$E_2 \leq (cw(\theta_n))^{-1} \sum_{i=1}^{N(K)} D(\rho)[\|\Delta_i^K\| + \|f_i^K\|] n^{-\rho/2} + 2(cw(\theta_n))^{-1} 2^{-K}.$$

By taking $b = b_n$ and $c = b_n h(\theta_n)/4$ and substituting (4.36) and (4.37) into (4.23), we have

(4.38)
$$\mathbb{P}(r_n|\hat{\theta}_n - \theta_n| > 2^M)$$
$$\leq C(\theta, \eta) D(\rho) b_n^{-1} \left(\frac{r_n}{n}\right)^{\rho/2} N(K) \sup_{f \in \mathcal{F}} \|f\| \sum_{j=M}^{J} 2^{(-1/2)j\rho}$$
$$+ C(\theta, \eta) D(\rho) N(K) b_n^{-1} n^{-\rho/2} \sup_{f \in \mathcal{F}} \|f\|$$
$$+ C(\theta, \eta) b_n^{-1} 2^{-K} + \mathbb{P}(|\hat{\theta}_n - \theta_n| > \eta) + \mathbb{P}(N(\delta_n) \leq b_n).$$

Choosing $K$ such that $b_n^{-1} 2^{-K} \sim \varepsilon_n$ and $r_n$ such that $N(K) b_n^{-1} r_n^{\rho/2} n^{-\rho/2} = 1$, we obtain

(4.39) $$\lim_{n \to +\infty} \mathbb{P}(r_n|\hat{\theta}_n - \theta_n| > 2^M) \leq C(\theta, \eta) D(\rho) \sup_{f \in \mathcal{F}} \|f\| \sum_{j=M}^{+\infty} 2^{(-1/2)j\rho}.$$

Finally, letting $M$ tend to infinity ends the proof.

REMARK 4. If $b_n \geq b > 0$, we choose $r_n = n$. Firstly, let $n$ go to infinity, then let $M$ go to infinity and finally let $K$ go to infinity to obtain $\lim_{M \to +\infty} \lim_{n \to +\infty} \mathbb{P}(n|\hat{\theta}_n - \theta_n| > 2^M) = 0$, without posing any restriction on the covering numbers.

**Acknowledgments.** We gratefully acknowledge the constructive comments of an Associate Editor and two anonymous referees that were very helpful for the improvement of the paper.

S. Ben Hariz  
Laboratoire de Statistique et Processus  
Département de Mathématiques  
Université du Maine  
Avenue Olivier Messiaen  
72085 Le Mans Cedex 9  
France  
E-mail: samir.ben_hariz@univ-lemans.fr

J. J. Wylie  
Q. Zhang  
Department of Mathematics  
City University of Hong Kong  
83 Tat Chee Avenue  
Kowloon Tong  
Hong Kong  
E-mail: wylie@math.cityu.edu.hk  
        mazq@cityu.edu.hk